\documentclass[12pt]{article}%
\usepackage{amsfonts}
\usepackage{amssymb}
\usepackage{graphicx}
\usepackage{amsmath}
\usepackage{makeidx}%
\providecommand{\U}[1]{\protect\rule{.1in}{.1in}}
\newtheorem{theorem}{Theorem}

\newtheorem{algorithm}[theorem]{Algorithm}

\newtheorem{conjecture}[theorem]{Conjecture}
\newtheorem{corollary}[theorem]{Corollary}

\newtheorem{problem}[theorem]{Problem}

\numberwithin{equation}{section}
\numberwithin{theorem}{section}
\makeindex

\ifx\pdfoutput\relax\let\pdfoutput=\undefined\fi
\newcount\msipdfoutput
\ifx\pdfoutput\undefined\else
\ifcase\pdfoutput\else
\ifx\paperwidth\undefined\else
\ifdim\paperheight=0pt\relax\else\pdfpageheight\paperheight\fi
\ifdim\paperwidth=0pt\relax\else\pdfpagewidth\paperwidth\fi
\fi\fi\fi
\begin{document}

\title{\textbf{Algorithms and Convergence Results of Projection Methods for
Inconsistent Feasibility Problems: A Review}}
\author{Yair Censor$^{1}$ and Maroun Zaknoon$^{2}\bigskip$\\$^{1}$Department of Mathematics, University of Haifa\\Mt.\ Carmel, Haifa 3498838, Israel\\\{yair@math.haifa.ac.il\}$\bigskip$\\$^{2}$Department of Mathematics\\The Arab Academic College for Education\\22 HaHashmal Street, Haifa 32623, Israel\\\{zaknoon@arabcol.ac.il\}}
\date{February 21, 2018. Revised: March 19, 2018. Revised: April 25, 2018\bigskip.\\
\textbf{We dedicate this paper to Adi Ben-Israel, our scientific father and
grandfather, respectively.}}
\maketitle

\begin{abstract}
The convex feasibility problem (CFP) is to find a feasible point in the
intersection of finitely many convex and closed sets. If the intersection is
empty then the CFP is inconsistent and a feasible point does not exist.
However, algorithmic research of inconsistent CFPs exists and is mainly
focused on two directions. One is oriented toward defining solution concepts
other that will apply, such as proximity function minimization wherein a
proximity function measures in some way the total violation of all
constraints. The second direction investigates the behavior of algorithms that
are designed to solve a consistent CFP when applied to inconsistent problems.
This direction is fueled by situations wherein one lacks a priory information
about the consistency or inconsistency of the CFP or does not wish to invest
computational resources to get hold of such knowledge prior to running his
algorithm. In this paper we bring under one roof and telegraphically review
some recent works on inconsistent CFPs.

\end{abstract}

\section{Introduction\label{introduction}}

\textbf{Inconsistent feasibility problems. }Feasibility problems require to
find a point in a given set $C$, any point, not a particular point such as,
for example, one that optimizes some given function over $C$, which would
constitute a problem of constrained optimization. Often times the set $C$ is
given as an intersection $C:=\cap_{i=1}^{m}C_{i}$ of a finite family of sets
$\{C_{i}\}_{i=1}^{m}.$ The convex feasibility problem (CFP) is to find a
feasible point $x^{\ast}\in C=\cap_{i=1}^{m}C_{i}$ when all sets $C_{i}$ are
convex and commonly also assumed to be closed. This prototypical problem
underlies the modeling of real-world problems in the set theoretic estimation
approach of Combettes%
\index{Combettes}
\cite{Combettes1993} such as convex set theoretic image recovery
\cite{Combettes1997}\ and many other fields, see, e.g., the pointers and
references in Bauschke%
\index{Bauschke}
and Borwein%
\index{Borwein}
\cite[Section 1]{Bauschke96} and in Cegielski's%
\index{Cegielski}
book \cite[Section 1.3]{Cegielski2012Book}. In this approach, constraints of
the real-world problem are represented by the demand that a solution should
belong to sets $C_{i},$ called constraint sets.

If $C\neq\emptyset$ does not hold\ then the CFP is inconsistent and a feasible
point does not exist. However, algorithmic research of inconsistent CFPs
exists and is mainly focused on two directions. One is oriented toward
defining solution concepts other than $x^{\ast}\in C=\cap_{i=1}^{m}C_{i}$ that
will apply, such as proximity function minimization wherein a proximity
function measures in some way the total violation of all constraints. The
second direction investigates the behavior of algorithms that are designed to
solve a consistent CFP when applied to inconsistent problems. The latter
direction is fueled by situations wherein one lacks a priory information about
the consistency or inconsistency of the CFP or does not wish to invest
computational resources to get hold of such knowledge prior to running his
algorithm. The next paragraphs on projection methods are quoted from the
introduction of Censor%
\index{Censor}
and Cegielski%
\index{Cegielski}
\cite{censor2015projection}.

\textbf{Projection methods}. Projections onto sets are used in a wide variety
of methods in optimization theory but not every method that uses projections
really belongs to the class of projection methods as we mean it here. Here
\textit{projection methods} are iterative algorithms that use projections onto
sets while relying on the general principle that when a family of (usually
closed and convex) sets is present then projections (or approximate
projections) onto the given individual sets are easier to perform than
projections onto other sets (intersections, image sets under some
transformation, etc.) that are derived from the given family of individual sets.

A projection algorithm reaches its goal, related to the whole family of sets,
by performing projections onto the individual sets. Projection
algorithms\textit{\ }employ projections (or approximate projections) onto
convex sets in various ways. They may use different kinds of projections and,
sometimes, even use different projections within the same algorithm. They
serve to solve a variety of problems which are either of the feasibility or
the optimization types. They have different algorithmic structures, of which
some are particularly suitable for parallel computing, and they demonstrate
nice convergence properties and/or good initial behavior patterns in some
significant fields of applications.

Apart from theoretical interest, the main advantage of projection methods,
which makes them successful in real-world applications, is computational. They
commonly have the ability to handle huge-size problems of dimensions beyond
which other, more sophisticated currently available, methods start to stutter
or cease to be efficient. This is so because the building bricks of a
projection algorithm are the projections onto the given individual sets
(assumed and actually easy to perform) and the algorithmic structures are
either sequential or simultaneous or in-between, such as in the
block-iterative projection (BIP) methods or in the more recent
string-averaging projection (SAP) methods. An advantage of projection methods
is that they work with initial data and do not require transformation of, or
other operations on, the sets describing the problem.

\textbf{Purpose of the paper}. We present an effort to bring under one roof
and telegraphically review some recent works on inconsistent CFPs. This should
be helpful to researchers, veterans or newcomers, by directing them to some of
the existing resources. The vast amount of research papers in the field of
projection methods\ makes it sometimes difficult to master even within a
specific sub-area. On the other hand, projection methods send branches both
into fields of applications wherein real-world problems are solved and into
theoretical areas in mathematics such as, but not only, fixed point theory and
variational inequalities. Researchers in each of these, seemingly
perpendicular, directions might benefit from this review.

\textbf{A word about notations}. We entertained the thought to unify all
notations but quickly understood that the game is not worth the
candle\footnote{Meaning that what we would get from this undertaking is not
worth the effort we would have to put into it.}. With notations left as they
appear in the original publications it will make it easier for a reader when
choosing to consult the original papers.

\textbf{An apology}. Oversight and lack of knowledge are human traits which we
are not innocent of. Therefore, we apologize for omissions and other
negligence and lacunas in this paper. We kindly ask our readers to communicate
to us any additional items and informations that fit the structure and spirit
of the paper and we will gladly consider those for inclusion in future
revisions, extensions and updates of the paper that we will post on arXiv.

\textbf{Organization of the paper}. Section \ref{sec:two} contains our review
divided into 17 subsections. Each subsection is focused on, and is centered
around, one or two historical or recent works. We use these \textquotedblleft
lead\textquotedblright\ references to organize the subsections chronologically
from older to recent works. An author index at the end of the paper will help
locate the results reviewed here.

\section{\textbf{Algorithms and Convergence results of Projection Methods for
Inconsistent Feasibility Problems}\label{sec:two}}

\subsection{1959: Composition of projections onto two disjoint convex
sets\label{subsec:CG}}

Cheney%
\index{Cheney}
and Goldstein%
\index{Goldstein}
\cite{Cheney1959} showed that if $K_{1}$ and $K_{2}$ are two closed and convex
subsets\ of a Hilbert space, and $P_{i}$ are the corresponding orthogonal
projections onto $K_{i},$ where $i=1,2$, then every fixed point of the
composition $Q:=P_{1}P_{2}$ is a point of $K_{1}$ closest to $K_{2}$.
Moreover, they showed that if one of the sets is compact or if one of the sets
is finite-dimensional and the distance is attained then a fixed point of $Q$
will be obtained by iterations of $Q$. In particular, if both sets are
polytopes in a finite-dimensional Euclidean space, the distance between two
sets is attained, and consequently a fixed point of $Q$ will be obtained by
iterations of $Q$. Their results are in the following three theorems.

\begin{theorem}
\label{Cheney1959 mark 1}\cite[Theorem 2]{Cheney1959} Let $K_{1}$ and $K_{2}$
be two closed convex sets in Hilbert space. Let $P_{i}$ denote the proximity
map for $K_{i}$. Any fixed point of $P_{1}P_{2}$ is a point of $K_{1}$ nearest
$K_{2}$, and conversely.
\end{theorem}

\begin{theorem}
\label{Cheney1959 mark 2}\cite[Theorem 4]{Cheney1959} Let $K_{1}$ and $K_{2}$
be two closed convex sets in Hilbert space and $Q$ the composition $P_{1}%
P_{2}$ of their proximity maps. Convergence of $Q^{n}x$ to a fixed point of
$Q$ is assured when either (a) one set is compact, or (b) one set is
finite-dimensional and the distance between the sets is attained.
\end{theorem}

\begin{theorem}
\label{Cheney1959 mark 3}\cite[Theorem 5]{Cheney1959} In a finite-dimensional
Euclidean space, the distance between two polytopes is attained, a polytope
being the intersection of a finite family of half-spaces.
\end{theorem}

In this connection, see also Theorems 4.1 and 4.2 in the paper by Kopeck\'{a}%
\index{Kopeck\'{a}}
and Reich%
\index{Reich}
\cite{kr2004}.

\subsection{1967: Cyclic convergence of sequential projections onto $m$ sets
with empty intersection}

Gubin%
\index{Gubin}%
, Polyak%
\index{Polyak}
and Raik%
\index{Raik}
\cite{Gubin1967} consider $m$ closed convex subsets, $C_{1},C_{2},\cdots
,C_{m}$, of\ a normed space $E$, and studied the behavior of the sequence
generated according the rule%
\begin{equation}
x^{n+1}=P_{i\left(  n\right)  }x^{n},\text{ where }i\left(  n\right)
:=n\left(  \operatorname{mod}m\right)  +1\label{Gubin1967 disp 1}%
\end{equation}
with $x^{0}$ arbitrary. They showed that if one of the sets is bounded, the
subsequences%
\begin{equation}
\left\{  x^{mn+1}\right\}  _{n\in%
%TCIMACRO{\U{2115} }%
%BeginExpansion
\mathbb{N}
%EndExpansion
},\left\{  x^{mn+2}\right\}  _{n\in%
%TCIMACRO{\U{2115} }%
%BeginExpansion
\mathbb{N}
%EndExpansion
},\cdots,\left\{  x^{mn+m}\right\}  _{n\in%
%TCIMACRO{\U{2115} }%
%BeginExpansion
\mathbb{N}
%EndExpansion
}%
\end{equation}
converge weakly to cluster points $\overline{x}_{1},\overline{x}_{2}%
,\cdots,\overline{x}_{m},$ respectively, that constitute a cycle, i.e.,%
\begin{equation}
\overline{x}_{2}=P_{2}\overline{x}_{1},\text{ }\overline{x}_{3}=P_{3}%
\overline{x}_{2},\text{ }\cdots\text{\ },\overline{x}_{m}=P_{m}\overline
{x}_{m-1},\text{ }\overline{x}_{1}=P_{1}\overline{x}_{m}%
.\label{Gubin1967 disp 2}%
\end{equation}

They proved the following (slightly paraphrased here) version of \cite[Theorem
2]{Gubin1967}.

\begin{theorem}
\label{Gubin1967 mark 1} Let all $C_{i},$ $i=1,2,\cdots,m,$ be closed, convex
and nonempty subsets of $E$ and at least one of them (for explicitness,
$C_{1}$) be bounded. Then it is possible to find points $\overline{x}_{i}\in
C_{i},$ $i=1,2,\cdots,m,$ such that $P_{i+1}\left(  \overline{x}_{i}\right)
=\overline{x}_{i+1},$ $i=1,2,\cdots,m-1,$ $P_{1}\left(  \overline{x}%
_{m}\right)  =\overline{x}_{1},$ while in the method ($\ref{Gubin1967 disp 1}%
$) we have $x^{km+i+1}-x^{km+i}\rightarrow\overline{x}_{i+1}-\overline{x}%
_{i},$ and $x^{km+i}$ weakly converges to $\overline{x}_{i}$ as $k\rightarrow
\infty.$ If, in addition, any of the following conditions is satisfied,
\newline(a) all $C_{i}$ with the possible exception of one ($C_{\overline{i}}%
$), are uniformly convex with the common function $\delta\left(  \tau\right)
;$\newline(b) $E$ is finite-dimensional;\newline(c) all $C_{i}$ are (closed)
half-spaces;\newline then the convergence will be strong. If all $C_{i},$
apart from possibly one, are also strongly convex, and $\cap_{i=1}^{m}C_{i}$
is empty, the sequence $x^{km+i}$ converges to $\overline{x}_{i}$ at the rate
of a geometrical progression.
\end{theorem}

\subsection{1983: The limits of the cyclic subsequences approach a single
point as relaxation goes to zero\label{susec:ceg83}}

Censor%
\index{Censor}%
, Eggermont%
\index{Eggermont}
and Gordon%
\index{Gordon}
\cite{Censor1983} investigate the behavior of Kaczmarz's%
\index{Kaczmarz}
method with relaxation for inconsistent systems. They show that when the
relaxation parameter goes to zero, the limits of the cyclic subsequences (See
Theorem \ref{Gubin1967 mark 1} here) generated by the method approach a
weighted least squares solution of the system. This point minimizes the sum of
the squares of the Euclidean distances to the hyperplanes of the system. If
the starting point is chosen properly, then the limits approach the minimum
norm weighted least squares solution. The proof is given for a block-Kaczmarz%
\index{Kaczmarz}
method.

Consider the linear system of equations $Ax=b,$ where $A\in%
%TCIMACRO{\U{211d} }%
%BeginExpansion
\mathbb{R}
%EndExpansion
^{m\times n}$, $b\in$ $%
%TCIMACRO{\U{211d} }%
%BeginExpansion
\mathbb{R}
%EndExpansion
^{m},$ $a_{i}$ is the $i$'th row of the matrix $A,$ and $b_{i}$ is the $i$'th
component the column vector $b.$ Kaczmarz's%
\index{Kaczmarz}
algorithm \cite{Kaczmarz1937} employs the iterative process%
\begin{equation}
x^{k+1}=x^{k}+\frac{b_{i}-\left\langle a_{i,}x^{k}\right\rangle }{\left\Vert
a_{i}\right\Vert ^{2}}a_{i}, \label{Censor1983 disp 1.2 apsence relaxation}%
\end{equation}
where $x^{0}\in%
%TCIMACRO{\U{211d} }%
%BeginExpansion
\mathbb{R}
%EndExpansion
^{n}$ is arbitrary, for solving the system $Ax=b$. Eggermont%
\index{Eggermont}%
, Herman%
\index{Herman}%
, and Lent%
\index{Lent}
\cite{Eggermont1981} rewrote the above $m$ and $n$ as $m=LM,$ $n=N,$with any
natural numbers $L,M,N$, partitioned $A$ and $b$ as%
\begin{equation}
A=\left(
\begin{array}
[c]{c}%
A_{1}\\
A_{2}\\
\vdots\\
A_{M}%
\end{array}
\right)  ,\ \ b=\left(
\begin{array}
[c]{c}%
b_{1}\\
b_{2}\\
\vdots\\
b_{M}%
\end{array}
\right)  , \label{Censor1983 disp 2.2}%
\end{equation}
with $A_{i}\in%
%TCIMACRO{\U{211d} }%
%BeginExpansion
\mathbb{R}
%EndExpansion
^{L\times N}$ and $b_{i}\in%
%TCIMACRO{\U{211d} }%
%BeginExpansion
\mathbb{R}
%EndExpansion
^{L},$ and proposed the following block-Kaczmarz%
\index{Kaczmarz}
method for solving a linear system of the above form.%
\begin{equation}
x^{0}\in%
%TCIMACRO{\U{211d} }%
%BeginExpansion
\mathbb{R}
%EndExpansion
^{N}\text{\ is arbitrary, }x^{k+1}=x^{k}+\lambda A_{i}^{T}\left(  b_{i}%
-A_{i}x^{k}\right)  ,\text{ }i=k\left(  \operatorname{mod}M\right)  +1,
\label{Censor1983 disp 2.3}%
\end{equation}
with relaxation $\lambda\in\left(  0,2\right)  .$

The effect of strong underrelaxation on the limits of the cyclic subsequences
generated by the block-Kaczmarz%
\index{Kaczmarz}
algorithm (\ref{Censor1983 disp 2.3}), investigated in \cite{Censor1983} is
included in the following theorem.

\begin{theorem}
\cite[Theorem 1]{Censor1983} For all $\lambda$ small enough,%
\begin{equation}
x^{\ast}\left(  \lambda\right)  =\underset{k\rightarrow\infty}{\lim}x^{kM}
\label{Censor1983 disp 3.1}%
\end{equation}
exists, and%
\begin{equation}
\underset{\lambda\rightarrow0}{\lim}x^{\ast}\left(  \lambda\right)  =A^{\dag
}b+\left(  Id-A^{\dag}A\right)  x^{0}. \label{Censor1983 disp 3.2}%
\end{equation}
(Where $A^{\dag}$ is the Moore-Penrsose%
\index{Moore-Penrsose}
inverse of the matrix $A.$)
\end{theorem}

This applies to every subsequence $\left\{  x^{kM+\ell}\right\}  _{k\geq0},$
$\ell\in\left\{  0,1,2,\cdots,M-1\right\}  $. In Censor%
\index{Censor}%
, Eggermont%
\index{Eggermont}
and Gordon%
\index{Gordon}
\cite[Page 91]{Censor1983}, it is shown that algorithm
(\ref{Censor1983 disp 1.2 apsence relaxation}) is a special case of algorithm
(\ref{Censor1983 disp 2.3}), with $L=1,$ and that the equality
(\ref{Censor1983 disp 3.2}) means $\ $that $\lim_{\lambda\rightarrow0}x^{\ast
}\left(  \lambda\right)  $ is a least squares solution of the system $Ax=b$. A
relevant remark concerning the behavior of the Cimmino method in the
inconsistent case appears in the Remark on pages 286--287 of the 1983 paper by
Reich%
\index{Reich}
\cite{reich83}.

\subsection{1993 and 1994: Alternating projection algorithms for two sets}

Bauschke%
\index{Bauschke}
and Borwein%
\index{Borwein}
\cite{Bauschke1993} investigated the convergence of the von Neumann's
\index{von Neumann}
alternating projection method for two arbitrary closed convex nonempty subsets
$A,$ $B$ of a Hilbert space $H.$ Finding a point in $A\cap B,$ or if $A\cap B$
is empty a good substitute for it, is a basic problem in various areas of mathematics.

Defining the distance between two nonempty subsets $M,$ $N$ by $d\left(
M,N\right)  :=\inf\left\{  \left\Vert m-n\right\Vert \mid m\in M,\ n\in
N\right\}  $ and denoting $E:=\left\{  a\in A\mid d\left(  a,B\right)
=d\left(  A,B\right)  \right\}  $ and $F:=\{b\in B\mid d\left(  b,A\right)
=d\left(  B,A\right)  \},$ one notes that if $A\cap B\neq\varnothing$ then
$E=F=$ $A\cap B.$ The projection of any point $x$ onto a closed convex
nonempty subset $C$ is denoted by $P_{C}x.$ The von Neumann%
\index{von Neumann}
algorithm for finding a point in $A\cap B$ is as follows: Given a starting
point $x\in X\subseteq H$, define, for every integer $n\geq1,$ the terms of
the sequences $\left(  a_{n}\right)  ,$ $\left(  b_{n}\right)  $ by%
\begin{equation}
b_{0}:=x,\ \ \ a_{n}:=P_{A}b_{n-1},\text{ \ \ }b_{n}:=P_{B}a_{n}.
\end{equation}
von Neumann proved that both sequences converge to $P_{A\cap B}\left(
x\right)  $ in norm when $A,$ $B$ are closed subspaces.

Assuming \cite[Page 201]{Bauschke1993} that $A,$ $B$\ are closed affine
subspaces, say $A=a+K,$ $B=b+L$ for vectors $a,b\in X$ and closed subspaces
$K,$ $L.$ The angle between $K$ and $L$ is denoted by $\gamma\left(
K,L\right)  $ Bauschke and Borwein proved the following.

\begin{theorem}
\cite[Theorem 4.11]{Bauschke1993} If $K+L$ is closed, then the von Neumann%
\index{Neumann}
sequences converge linearly with rate $\cos\gamma\left(  K,L\right)  $
independent of the starting point. In particular, this happens whenever one of
the following conditions holds: (i) $K$ or $L$ has finite dimension, (ii) $K$
or $L$ has finite codimension.
\end{theorem}

In \cite{Bauschke1994} Bauschke%
\index{Bauschke}
and Borwein%
\index{Borwein}%
\ analyzed Dykstra's%
\index{Dykstra}
algorithm for two arbitrary closed convex sets in a Hilbert space $X$. They
greatly expanded on the Cheney%
\index{Cheney}
and Goldstein%
\index{Goldstein}
papers (infinite dimensions, characterizations, etc.). See Subsections
\ref{subsec:CG} here and the recent work of Kopeck\'{a}%
\index{Kopeck\'{a}}
and Reich%
\index{Reich}
\cite{Kopecka2012} in Subsection \ref{subsec:KR} here.

\subsection{1994: Least-squares solutions of inconsistent signal feasibility
problems in a product space}

Combettes's%
\index{Combettes}
\cite{Combettes1994} presents parallel projection methods to find
least-squares solutions to inconsistent convex set theoretic signal synthesis
problems. The problem of finding a signal that minimizes a weighted average of
the squares of the distances to constraint sets is reformulated in a product
space, where it is equivalent to that of finding a point that lies in a
particular subspace and at minimum distance from the Cartesian product of the
original sets. A solution is obtained in the product space via methods of
alternating projections which naturally lead to methods of parallel
projections in the original space. The convergence properties of the proposed
methods are analyzed and signal synthesis applications are demonstrated.

The, possibly inconsistent, feasibility problem: Find $a^{\ast}\in\cap
_{i=1}^{m}S_{i},$ where the $S_{i}$s are closed and convex subsets of a
Hilbert space $\Xi,$ is replaced by the unconstrained weighted least-squares
minimization problem%
\begin{equation}
\min\text{ }\{\text{ }\Phi\left(  a\right)  \mid a\in\Xi\},
\label{Combettes1994 disp 5}%
\end{equation}
where $\Phi\left(  a\right)  :=\frac{1}{2}\sum_{i=1}^{m}w_{i}d\left(
a,S_{i}\right)  ^{2},$ $d\left(  a,S_{i}\right)  :=\inf\left\{  d\left(
a,b\right)  \mid b\in S_{i}\right\}  ,$ and $\left(  w_{i}\right)  _{1\leq
i\leq m}$ are strictly convex weights, i.e., $\sum_{i=1}^{m}w_{i}=1$ and
$\forall i\in\left\{  1,\cdots,m\right\}  $ $\omega_{i}>0.$ In other words,
the goal is to solve%
\begin{equation}
\text{Find }a^{\ast}\in G:=\left\{  a\in\Xi\mid\text{ }\Phi\left(  a\right)
\leq\Phi\left(  b\right)  \text{ for all }b\in\Xi\right\}  .
\label{Combettes1994 disp 7}%
\end{equation}
In the Cartesian product space $\Xi^{m},$ with the scalar product
$\left\langle \left\langle \mathbf{a},\mathbf{b}\right\rangle \right\rangle
:=\sum_{i=1}^{m}w_{i}\left\langle a^{\left(  i\right)  },b^{\left(  i\right)
}\right\rangle $ for all $\mathbf{a:=}\left(  a^{\left(  1\right)
},a^{\left(  2\right)  },\cdots,a^{\left(  m\right)  }\right)  \in\Xi^{m}$ and
$\mathbf{b:=}(b^{\left(  1\right)  },$ $b^{\left(  2\right)  },$ $\cdots,$
$b^{\left(  m\right)  })\in\Xi^{m}$ the problem (\ref{Combettes1994 disp 7})
is reformulated as%
\begin{equation}
\text{Find }\mathbf{a}^{\ast}\in\mathbf{G:=}\left\{  \mathbf{a\in D\mid
}d\left(  \mathbf{a},\mathbf{S}\right)  =d\left(  \mathbf{D},\mathbf{S}%
\right)  \right\}  =\left\{  \mathbf{a\in D\mid}P_{\mathbf{D}}\left(
P_{\mathbf{S}}\left(  \mathbf{a}\right)  \right)  =\mathbf{a}\right\}  ,
\end{equation}
where $\mathbf{D=}\left\{  \left(  a,a,\cdots,a\right)  \in\Xi^{m}\mid a\in
\Xi\right\}  $, $\mathbf{S=}S_{1}\times S_{2}\cdots\times S_{m}$, and
$P_{\mathbf{D}},$ $P_{\mathbf{S}}$ are the orthogonal projections onto the
sets $\mathbf{D}$ and $\mathbf{S,}$ respectively. All quantities related to
the product space are written in boldface symbols. Solving this problem using
two methods for finding a fixed point of the composition $P_{\mathbf{D}}\circ
P_{\mathbf{S}}$ and translating back the results to the original space $\Xi,$
the following two convergence results are obtained.

\begin{theorem}
\label{Combettes1994 mark 1 Theo 4}\cite[Theorem 4]{Combettes1994} Suppose
that one of the $S_{i}$s is bounded. Then, for any $a_{0}\in\Xi,$ every
sequence of iterates $\left(  a_{n}\right)  _{n\geq0}$ defined by%
\begin{equation}
\ a_{n+1}=a_{n}+\lambda_{n}\left(  \sum_{i=1}^{m}w_{i}P_{i}\left(
a_{n}\right)  -a_{n}\right)  ,\label{Combettes1994 disp 1 Eq 32}%
\end{equation}
where $\left(  \lambda_{n}\right)  _{n\geq0}\subseteq\left[  \varepsilon
,2-\varepsilon\right]  $ with $0<\varepsilon<1,$ converges weakly to a point
in $G.$
\end{theorem}

Assuming that $\left(  \alpha_{n}\right)  _{n\geq0}$ fulfills%
\begin{equation}
\lim_{n\rightarrow+\infty}\alpha_{n}=1,\text{ }%
%TCIMACRO{\tsum \limits_{n\geq0}}%
%BeginExpansion
{\textstyle\sum\limits_{n\geq0}}
%EndExpansion
\left(  1-\alpha_{n}\right)  =+\infty\text{ and}\lim_{n\rightarrow+\infty
}\left(  \alpha_{n+1}-\alpha_{n}\right)  \left(  1-\alpha_{n+1}\right)
^{-2}=0 \label{Combettes1994 disp 2 Eq 12}%
\end{equation}
the next theorem holds.

\begin{theorem}
\label{Combettes1994 mark 2 Theo 5}\cite[Theorem 5]{Combettes1994} Suppose
that one of the $S_{i}$s is bounded. Then, for any $a_{0}\in\Xi,$ every
sequence of iterates $\left(  a_{n}\right)  _{n\geq0}$ defined by%
\begin{equation}
\ a_{n+1}=\left(  1-\alpha_{n}\right)  a_{0}+\alpha_{n}\left(  \lambda
\sum_{i=1}^{m}w_{i}P_{i}\left(  a_{n}\right)  -\left(  1-\lambda\right)
a_{n}\right)  ,\label{Combettes1994 disp 3 Eq 37}%
\end{equation}
where $\left(  \alpha_{n}\right)  _{n\geq0}$ is as in
(\ref{Combettes1994 disp 2 Eq 12}) and $0<\lambda\leq2,$ converges strongly to
the projection of $a_{0}$ onto $G.$
\end{theorem}

\subsection{1995 and 2003: The method of cyclic projections for closed convex
sets in Hilbert space}

Bauschke%
\index{Bauschke}%
, Borwein%
\index{Borwein}
and Lewis%
\index{Lewis}
\cite{BauschkeBS1997} consider closed convex nonempty sets $C_{1},C_{2}%
,\cdots,C_{N}$ in a real Hilbert space $H,$ with corresponding projections
$P_{1},P_{2},\cdots,P_{N}$ and systematically study composition of projections
in the inconsistent case. For an arbitrary starting point $x^{0}\in H$ the
method of cyclic projections generates $N$ sequences $\left(  x_{i}%
^{n}\right)  _{n}$ by%
\begin{equation}%
\begin{array}
[c]{l}%
x_{1}^{1}:=P_{1}x^{0},\text{ }x_{2}^{1}:=P_{2}x_{1}^{1},\text{ }\cdots,\text{
}x_{N}^{1}:=P_{N}x_{N-1}^{1},\\
x_{1}^{2}:=P_{1}x_{N}^{1},\text{ }x_{2}^{2}:=P_{2}x_{1}^{2},\text{ }%
\cdots,\text{ }x_{N}^{2}:=P_{N}x_{N-1}^{2},\\
x_{1}^{3}:=P_{1}x_{N}^{2},\text{ \ \ \ \ \ \ \ \ \ }\cdots\cdots\text{ .}%
\end{array}
\end{equation}

They collected these sequences cyclically in one sequence $(x^{0},$ $x_{1}%
^{1},$ $x_{2}^{1},$ $\cdots,$ $x_{N}^{1},$ $x_{1}^{2},$ $x_{2}^{2},$ $\cdots)$
to which they referred as the orbit generated by $x^{0}$ or the orbit with
starting point $x^{0}.$ They further defined the composite projections
operators%
\begin{equation}
Q_{1}:=P_{1}P_{N}P_{N-1}\cdots P_{2},\text{ }Q_{2}:=P_{2}P_{1}P_{N}\cdots
P_{3},\text{ }\cdots\text{ },\text{ }Q_{N}:=P_{N}P_{N-1}\cdots P_{1}%
\end{equation}
which allows to write more concisely%
\begin{equation}
x_{i}^{n}:=Q_{i}^{n-1}x_{i}^{1},\text{ for all }n\geq1\text{ and for every }i;
\end{equation}
after setting $P_{0}:=P_{N},$ $P_{N+1}:=P_{1},$ $x_{0}^{n}:=x_{N}^{n-1},$ and
$x_{N+1}^{n}:=x_{1}^{n+1},$ they reached%
\begin{equation}
x_{i+1}^{n}=P_{i+1}x_{i}^{n},\text{ for all }n\geq1\text{ and every }i.
\end{equation}
When appropriate, they similarly identify $i=0$ with $i=N$ and $i=N+1$ with
$i=1.$

They gave a dichotomy result on orbits which roughly says that if each $Q_{i}$
is fixed point free then, the orbit has no bounded subsequence; otherwise,
each subsequence $\left(  x_{i}^{n}\right)  $ converges weakly to some fixed
point of $Q_{i}.$ Two central questions were posed:

\begin{enumerate}
\item When does each $Q_{i}$ have a fixed point?

\item If each $Q_{i}$ has a\ fixed point, when do the subsequences $\left(
x_{i}^{n}\right)  $ converge in norm (or even linearly)?
\end{enumerate}

Concerning Question 1, They provide sufficient conditions for the existence of
fixed points or approximate fixed points (that is $\inf_{x\in H}\left\Vert
x-Q_{i}x\right\Vert =0,$ for each $i$). It follows that while fixed points of
$Q_{i}$ need not exist for non-intersecting closed affine subspaces,
approximate fixed points must.

In respect to Question 2, a variety of conditions guaranteeing norm
convergence (in the presence of fixed points for each $Q_{i}$) is offered: one
of the sets $C_{i}$ has to be (boundedly) compact or all sets are
\index{Kadec/Klee}
or convex polyhedra, or affine subspaces. In the affine subspace case each
sequence $\left(  x_{i}^{n}\right)  $ converges to the fixed point of $Q_{i}$
nearest to $x^{0}.$ Moreover the convergence is linear, whenever the angle of
the $N$-tuple of the associated closed subspaces is positive.

In subsequent work, Bauschke%
\index{Bauschke}
\cite{Bauschke2003} showed that the composition of finitely many projections
$P_{N},P_{N-1},\cdots,P_{1}$ is asymptotically regular, i.e.,%
\begin{equation}
\left(  P_{N},P_{N-1},\cdots,P_{1}\right)  ^{k}x-\left(  P_{N},P_{N-1}%
,\cdots,P_{1}\right)  ^{k+1}x\rightarrow0,\text{ for every }x\in X.
\end{equation}
thus proving the so-called \textquotedblleft zero displacement
conjecture\textquotedblright\ of Bauschke%
\index{Bauschke}%
, Borwein%
\index{Borwein}
and Lewis%
\index{Lewis}
\cite{BauschkeBS1997}.

\subsection{1999: Hard-constrained inconsistent signal feasibility
problems\label{subsec:hard-soft}}

Combettes%
\index{Combettes}
and Bondon%
\index{Bondon}
\cite{Combettes1999} consider the problem of synthesizing feasible signals in
a Hilbert space $\mathcal{H}$, with inconsistent convex constraints that are
divided into two parts, the hard constraints and the soft constraints. They
look for a point in $\mathcal{H}$ which satisfies the hard constraints
imperatively and minimizes the violation of the soft constraints.

Denote by $\Gamma$ the class of all lower semicontinuous proper convex
functions from $\mathcal{H}$ into $\left]  -\infty,+\infty\right]  $. Given
$g\in\Gamma$ and $\alpha\in%
%TCIMACRO{\U{211d} }%
%BeginExpansion
\mathbb{R}
%EndExpansion
,$ the closed and convex set $lev_{\leq\alpha}g:=\left\{  x\in\mathcal{H}%
\ \mid g\left(  x\right)  \leq\alpha\right\}  $ is the lower level set of $g$
at height $\alpha,$ and the nonempty convex set $domg:=\left\{  x\in
\mathcal{H}\ \mid g\left(  x\right)  <+\infty\right\}  $ is its domain. The
goal of a convex set theoretic signal synthesis (design or estimation) problem
in $\mathcal{H}$ is to produce a signal $x^{\ast}$ that satisfies convex
constraints, say,%
\begin{equation}
\text{find }x^{\ast}\in S=\cap_{i\in I}S_{i}\text{, where }\left(  \forall
i\in I\right)  \text{ }S_{i}=lev_{\leq0}g_{i}, \label{Combettes1999 disp 1}%
\end{equation}
where $I$ is a finite index set, and $\left(  g_{i}\right)  _{i\in I}%
\subset\Gamma.$

Let $I^{\blacktriangle}\subset I$ denote the, possibly empty, hard constraints
index set, $I^{\triangle}=I\smallsetminus I^{\blacktriangle}$ the nonempty
soft constraints index set, $S^{\blacktriangle}=\cap_{i\in I^{\blacktriangle}%
}S_{i}$ the hard feasibility set and, by convention, $S^{\blacktriangle
}=\mathcal{H}$ if $I^{\blacktriangle}=\varnothing$. $S^{\triangle}=\cap_{i\in
I^{\triangle}}S_{i},$ $D^{\triangle}=\cap_{i\in I^{\triangle}}domg_{i}$, and
assume that $S^{\blacktriangle}\cap D^{\triangle}\neq\varnothing.$
$\mathcal{F}$ is the class of all increasing convex functions from $\left[
0,+\infty\right[  $ into $\left[  0,+\infty\right[  $ that vanish (only) at
$0$; every $f\in\mathcal{F}$ is extended to the argument $+\infty$ by setting
$f\left(  +\infty\right)  =+\infty.$ For every $g\in\Gamma,$ $g^{+}%
=\max\left\{  0,g\right\}  .$ The amount of violation of the soft constraints
$\left(  g_{i}\left(  x\right)  \leq0\right)  _{i\in I^{\triangle}}$ is
measured by an objective function $\Phi^{\triangle}:\mathcal{H\rightarrow
}\left[  0,+\infty\right]  $ of the general form%
\begin{equation}
\Phi^{\triangle}:=\sum_{i\in I^{\triangle}}f_{i}\circ g_{i}^{+},\text{ where
}\left(  f_{i}\right)  _{i\in I^{\triangle}}\subset\mathcal{F}.
\label{Combettes1999 disp 2}%
\end{equation}

The hard-constrained signal feasibility problem is to minimize the objective
$\Phi^{\triangle}$ of (\ref{Combettes1999 disp 2}) over the hard feasibility
set $S^{\blacktriangle}.$ Setting $\alpha^{\ast}=\inf_{x\in S^{\blacktriangle
}}\Phi^{\triangle}(x),$ the problem reads%
\begin{equation}
\text{find }x^{\ast}\in G:=\left\{  x\in S^{\blacktriangle}\mid\Phi
^{\triangle}\left(  x\right)  =\alpha^{\ast}\right\}  .
\label{Combettes1999 disp 3}%
\end{equation}

Combettes%
\index{Combettes}
and Bondon%
\index{Bondon}
\cite{Combettes1999} supply convergence theorems for the following processes
under various conditions ($P^{\blacktriangle}$ denotes the projector onto
$S^{\blacktriangle}$):%
\begin{equation}
\left(  \forall n\in%
%TCIMACRO{\U{2115} }%
%BeginExpansion
\mathbb{N}
%EndExpansion
\right)  \ x_{n+1}=P^{\blacktriangle}\left(  x_{n}-\gamma\nabla\Phi
^{\triangle}\left(  x_{n}\right)  \right)  ;
\end{equation}%
\begin{equation}
\left(  \forall n\in%
%TCIMACRO{\U{2115} }%
%BeginExpansion
\mathbb{N}
%EndExpansion
\right)  \text{ \ }x_{n+1}=\left(  1-\lambda_{n}\right)  x_{n}+\lambda
_{n}P^{\blacktriangle}\left(  x_{n}-\gamma\nabla\Phi^{\triangle}\left(
x_{n}\right)  \right)  ;
\end{equation}
and%
\begin{equation}
\left(  \forall n\in%
%TCIMACRO{\U{2115} }%
%BeginExpansion
\mathbb{N}
%EndExpansion
\right)  \text{ \ }x_{n+1}=\left(  1-\lambda_{n}\right)  r+\lambda
_{n}P^{\blacktriangle}\left(  x_{n}-\gamma\nabla\Phi^{\triangle}\left(
x_{n}\right)  \right)  ,
\end{equation}
for a fixed given $r\in\mathcal{H}$.

Defining the proximity function $\Phi^{\triangle}:=\frac{1}{2}\sum_{i\in
I^{\triangle}}w_{i}d\left(  \cdot,S_{i}\right)  ^{2}$, with weights $\left(
w_{i}\right)  _{i\in I^{\triangle}}\subset\left]  0,1\right]  ,$ $\sum_{i\in
I^{\triangle}}w_{i}=1,$ and denoting by $P_{i}$ the projector onto $S_{i},$
the authors prove that, under certain conditions, the iterative process
\begin{equation}
\left(  \forall n\in%
%TCIMACRO{\U{2115} }%
%BeginExpansion
\mathbb{N}
%EndExpansion
\right)  \text{ \ }x_{n+1}=\left(  1-\lambda_{n}\right)  x_{n}+\lambda
_{n}P^{\blacktriangle}\left(  \left(  1-\gamma\right)  x_{n}+\gamma\sum_{i\in
I^{\triangle}}w_{i}P_{i}\left(  x_{n}\right)  \right)
\end{equation}
generates sequences that converge, weakly or strongly, to a solution of
the\ hard-constrained signal feasibility problem, i.e., to a point in $G$ of
(\ref{Combettes1999 disp 3}) above. Also, they show that the iterative
process
\begin{equation}
\left(  \forall n\in%
%TCIMACRO{\U{2115} }%
%BeginExpansion
\mathbb{N}
%EndExpansion
\right)  \text{ \ }x_{n+1}=\frac{n}{n+1}P^{\blacktriangle}\left(  \left(
1-\gamma\right)  x_{n}+\gamma\sum_{i\in I^{\triangle}}w_{i}P_{i}\left(
x_{n}\right)  \right)
\end{equation}
generates sequences that converge strongly to $P_{G}\left(  0\right)  .$ In
all these iterative processes there appears $P^{\blacktriangle},$ the
projector onto $S^{\blacktriangle},$ which can potentially hinder practical
applications if this projection is not simple to calculate.

\subsection{2001: De Pierro's%
\index{De Pierro}
conjecture}

Bauschke%
\index{Bauschke}
and Edwards%
\index{Edwards}
\cite{BauschkeAndEdwards2005} describe De Pierro's%
\index{De Pierro}
conjecture as follows. Suppose we are given finitely many nonempty closed
convex sets in a real Hilbert space and their associated projections. For
suitable arrangements of the sets, it is known that the sequence obtained by
iterating the composition of the underrelaxed projections is weakly
convergent. The question arises how these weak limits vary as the
underrelaxation parameter tends to zero. In 2001, De Pierro%
\index{De Pierro}
conjectured \cite{DePierro2001} that the weak limits approach the least
squares solution nearest to the starting point of the sequence. In fact, the
result by Censor%
\index{Censor}%
, Eggermont%
\index{Eggermont}%
, and Gordon%
\index{Gordon}
\cite{Censor1983} described here in Subsection \ref{susec:ceg83}, implies De
Pierro's%
\index{De Pierro}
conjecture for affine subspaces in Euclidean space.

De Pierro's%
\index{De Pierro}
conjecture \cite[Section 3, Conjecture II]{DePierro2001} is succinctly
formulated in Bauschke%
\index{Bauschke}
and Edwards%
\index{Edwards}
\cite[Conjecture 1.6]{BauschkeAndEdwards2005} as follows. For a convex
feasibility problem with $N$ sets $(C_{i})_{i=1}^{N}$ in a Hilbert space $X,$
define, fort every $\lambda\in]0,1],$ the composition of underrelaxed
projections%
\begin{equation}
Q_{\lambda}:=((1-\lambda)Id+\lambda P_{C_{N}})\cdots((1-\lambda)Id+\lambda
P_{C_{2}})((1-\lambda)Id+\lambda P_{C_{1}}),
\end{equation}
where $Id$ is the identity and $P_{C_{i}}$ is the projection onto $C_{i}.$ The
corresponding sets of fixed points are defined by%
\begin{equation}
F_{\lambda}:=\operatorname*{Fix}Q_{\lambda}:=\{x\in X\mid x=Q_{\lambda}(x)\}
\end{equation}
and the aim is to understand the behavior of the sequence $(Q_{\lambda}%
^{n}(x))_{n\in%
%TCIMACRO{\U{2115} }%
%BeginExpansion
\mathbb{N}
%EndExpansion
}$ in terms of $\lambda\in]0,1]$, for an arbitrary $x\in X$.

\begin{conjecture}
(De Pierro%
\index{De Pierro}%
) \cite[Conjecture 1.6]{BauschkeAndEdwards2005}. Suppose that $F_{\lambda}%
\neq\emptyset$ for every $\lambda\in]0,1].$ Denoting, for all $x\in X$ and all
$\lambda\in]0,1]$ the limits $x_{\lambda}=\operatorname*{weak}\lim
_{n\rightarrow+\infty}Q_{\lambda}^{n}(x),$ De Pierro%
\index{De Pierro}
conjectured that $\lim_{\lambda\rightarrow0^{+}}x_{\lambda}=P_{\mathcal{L}%
}(x),$ where $\mathcal{L}$ is the set of least squares solutions of the convex
feasibility problem, i.e.,%
\begin{equation}
\mathcal{L}:=\left\{  x\in X\mid%
%TCIMACRO{\dsum \limits_{i=1}^{N}}%
%BeginExpansion
{\displaystyle\sum\limits_{i=1}^{N}}
%EndExpansion
\left\Vert x-P_{C_{i}}(x)\right\Vert ^{2}=\inf_{y\in X}%
%TCIMACRO{\dsum \limits_{i=1}^{N}}%
%BeginExpansion
{\displaystyle\sum\limits_{i=1}^{N}}
%EndExpansion
\left\Vert y-P_{C_{i}}(y)\right\Vert ^{2}\right\}  .
\end{equation}

\end{conjecture}

Bauschke%
\index{Bauschke}
and Edwards%
\index{Edwards}
proved this conjecture for families of closed affine subspaces satisfying a
metric regularity condition. Baillon%
\index{Baillon}%
, Combettes%
\index{Combettes}
and Cominetti%
\index{Cominetti}
\cite[Theorem 3.3]{Baillon2014} proved the conjecture under a mild geometrical
condition. Recently, Cominetti%
\index{Cominetti}%
, Roshchina%
\index{Roshchina}
and Williamson%
\index{Williamson}
\cite[Theorem 1]{CominettiTechRep2018} proved that this conjecture is false in
general by constructing a system of three compact convex sets in $%
%TCIMACRO{\U{211d} }%
%BeginExpansion
\mathbb{R}
%EndExpansion
^{3}$ for which the least squares solution exists but the conjecture fails to hold.

\subsection{2001: Proximity function minimization using multiple Bregman%
\index{Bregman}
projections\label{subsec:prox2001}}

Motivated by the geometric alternating minimization approach of Csisz\'{a}r%
\index{Csisz\'{a}r}
and Tusn\'{a}dy%
\index{Tusn\'{a}dy}
\cite{Csiszar1984} and the product space formulation of Pierra%
\index{Pierra}
\cite{Pierra1984}, Byrne%
\index{Byrne}
and Censor%
\index{Censor}
\cite{Byrne2001} derive a new simultaneous multiprojection algorithm that
employs generalized projections of Bregman%
\index{Bregman}
\cite{Bregman1967}, see also, e.g., Censor%
\index{Censor}
and Lent%
\index{Lent}
\cite{Censor1981} and Bauschke%
\index{Bauschke}
and Borwein%
\index{Borwein}
\cite{Bauschke1997}, to solve the convex feasibility problem (CFP) or, in the
inconsistent case, to minimize a proximity function that measures the average
distance from a point to all convex sets. For background material on Bregman%
\index{Bregman}
functions and Bregman%
\index{Bregman}
distances and projections see, e.g., the book of Censor and Zenios%
\index{Zenios}
\cite{censor1997Book}, \cite{Bauschke1997}, Solodov%
\index{Solodov}
and Svaiter%
\index{Svaiter}
\cite{Solodov2000}, Eckstein%
\index{Eckstein}
\cite{Eckstein1998}, \cite{Eckstein2003}, to name but a few. Byrne%
\index{Byrne}
and Censor%
\index{Censor}
\cite{Byrne2001} assume that the Bregman%
\index{Bregman}
distances involved are jointly convex, so that the proximity function itself
is convex. When the intersection of the convex sets is empty, but the closure
$\mathrm{cl}F(x)$ of the proximity function $F(x),$ defined by $\mathrm{cl}%
F(x):=\mathrm{lim\,inf}_{y\rightarrow x}\,F(y),$ has a unique global
minimizer, the sequence of iterates converges to this unique minimizer.
Special cases of this algorithm include the \textquotedblleft Expectation
Maximization Maximum Likelihood\textquotedblright\ (EMML) method in emission
tomography and a new convergence result for an algorithm that solves the split
feasibility problem.

Let $C_{i},$ $i=1,2,\ldots,I,$ be closed convex sets in the $J$-dimensional
Euclidean space $%
%TCIMACRO{\U{211d} }%
%BeginExpansion
\mathbb{R}
%EndExpansion
^{J}$ and let $C$ be their intersection. Let $S$ be an open convex subset of $%
%TCIMACRO{\U{211d} }%
%BeginExpansion
\mathbb{R}
%EndExpansion
^{J}$ and $f$ a Bregman%
\index{Bregman}
function from the closure $\overline{S}$ of $S$ into $%
%TCIMACRO{\U{211d} }%
%BeginExpansion
\mathbb{R}
%EndExpansion
$; see, e.g., \cite[Chapter 2]{censor1997Book}. For a Bregman%
\index{Bregman}
function $f(x)$, the Bregman%
\index{Bregman}
distance $D_{f}$ is defined by%
\begin{equation}
D_{f}(z,x):=f(z)-f(x)-\langle\nabla f(x),z-x\rangle,
\label{Byrne2001 disp 1.1}%
\end{equation}
where $\nabla f(x)$ is the gradient of $f$ at $x$. If the function $f$ has the
form $f(x)=\sum_{j=1}^{J}g_{j}(x_{j})$, with the $g_{j}$ scalar Bregman%
\index{Bregman}
functions, then $f$ and the associated $D_{f}(z,x)$ are called
\textit{separable}. With $g_{j}(t)=t^{2}$, for all $j$, the function
$f(x)=\sum_{j=1}^{J}g_{j}(x_{j})=\sum_{j=1}^{J}x_{j}^{2}$ is a separable
Bregman%
\index{Bregman}
function and $D_{f}(z,x)$ is the squared Euclidean distance between $z$ and
$x$. For each $i$, denote by $P_{C_{i}}^{f}(x)$ the \textit{Bregman%
\index{Bregman}
projection} of $x\in S$ onto $C_{i}$ with respect to the \textit{Bregman%
\index{Bregman}
function} $f,$ i.e., $D_{f}(P_{C_{i}}^{f}(x),x)\leq D_{f}(z,x)$, for all $z\in
C_{i}\cap\overline{S}$. In \cite[Eq. (1.2)]{Byrne2001} the proximity function
$F(x)$ is of the form%
\begin{equation}
F(x)=\sum_{i=1}^{I}D_{f_{i}}(P_{C_{i}}^{f_{i}}(x),x),
\label{Byrne2001 disp 1.2}%
\end{equation}
where the $D_{f_{i}}$ are Bregman%
\index{Bregman}
distances derived from possibly distinct, possibly nonseparable Bregman%
\index{Bregman}
functions $f_{i}$ with zones $S_{f_{i}}.$ The function $F$ is defined, for all
$x$ in the open convex set $U:=\cap_{i=1}^{I}S_{f_{i}},$ which is assumed
nonempty. The proximity function $F(x)$ of (\ref{Byrne2001 disp 1.2}) is
extended to all of $%
%TCIMACRO{\U{211d} }%
%BeginExpansion
\mathbb{R}
%EndExpansion
^{J}$ by defining $F(x)=+\infty,$ for all $x\notin U$ and its closure
$\mathrm{cl}F$ is as defined above. They proved convergence of their iterative
method whenever $\mathrm{cl}F$ has a unique minimizer or when the set
$C\cap{\overline{U}}$ is nonempty. The following algorithm is proposed.

\begin{algorithm}
\label{Byrne2001 mark 1 alg 4.1}\cite[Algorithm 4.1]{Byrne2001}.\newline%
\textbf{Initialization:} $x^{0}\in U$ is arbitrary.\newline\textbf{Iterative
Step:} Given $x^{k}$ find, for all $i=1,2,\ldots,I$, the projections
$P_{C_{i}}^{f_{i}}(x^{k})$ and calculate $x^{k+1}$ from%
\begin{equation}
\sum_{i=1}^{I}\nabla^{2}f_{i}(x^{k+1})x^{k+1}=\sum_{i=1}^{I}\nabla^{2}%
f_{i}(x^{k+1})P_{C_{i}}^{f_{i}}(x^{k}), \label{Byrne2001 disp 4.1}%
\end{equation}
where $\nabla^{2}f_{i}(x^{k+1})$ denotes the Hessian matrix (of second partial
derivatives) of the function $f_{i}$ at $x^{k+1}$.
\end{algorithm}

Let $F(x)$ be defined for $x\in U$ by (\ref{Byrne2001 disp 1.2}) and for other
$x\in%
%TCIMACRO{\U{211d} }%
%BeginExpansion
\mathbb{R}
%EndExpansion
^{J}$ let it be equal to $+\infty$ and let the set of minimizers of
$\mathrm{cl}F$ over $%
%TCIMACRO{\U{211d} }%
%BeginExpansion
\mathbb{R}
%EndExpansion
^{J}$ be denoted by $\Phi.$ Denote $\Gamma:=\mathrm{inf}\{\mathrm{cl}%
F(x)\,\mid x\in${$%
%TCIMACRO{\U{211d} }%
%BeginExpansion
\mathbb{R}
%EndExpansion
$}${^{J}}\}$ and consider the following assumptions.

\textbf{Assumption A1:} (\textit{Zone Consistency}) For every $i=1,2,\ldots
,I,$ if $x^{k}\in S_{f_{i}}$ then $P_{C_{i}}^{f_{i}}(x^{k})\in S_{f_{i}}.$

\textbf{Assumption A2:}\ For every $k=1,2,\ldots,$ the function $F_{k}%
(x):=\sum_{i=1}^{I}$ $D_{f_{i}}(P_{C_{i}}^{f_{i}}(x^{k}),x)$ has a unique
minimizer within $U.$

\textbf{Assumption A3:}\ If $\mathrm{cl}F(x)=0$ for some $x$ then $x$ is in
$C\cap{\bar{U}.}$

\begin{theorem}
\label{Byrne2001 mark 2 theo 4.1} \cite[Theorem 4.1]{Byrne2001} Let
Assumptions \textbf{A1}, \textbf{A2} and \textbf{A3} hold and assume that the
distances $D_{f_{i}}$ are jointly convex, for all $i=1,2,\cdots,I$, In
addition, assume that the set $\Phi$ is nonempty. If $\mathrm{cl}F$ has a
unique minimizer then any sequence $\left\{  x^{k}\right\}  ,$ generated by
Algorithm \ref{Byrne2001 mark 1 alg 4.1}, converges to this minimizer. If
$\Phi$ is not a singleton but $\Gamma=\mathrm{inf}\{\mathrm{cl}F(x)\,\mid
x\in{%
%TCIMACRO{\U{211d} }%
%BeginExpansion
\mathbb{R}
%EndExpansion
^{J}}\}=0$, then the intersection $C$ of the sets $C_{i}$ is nonempty and
$\left\{  x^{k}\right\}  $ converges to a solution of the CFP.
\end{theorem}

\subsection{2003: String-averaging projection schemes for inconsistent convex
feasibility problems}

Censor%
\index{Censor}
and Tom%
\index{Tom}
\cite{censor2003} study iterative projection algorithms for the convex
feasibility problem of finding a point in the intersection of finitely many
nonempty, closed and convex subsets in the Euclidean space. They propose
(without proof) an algorithmic scheme which generalizes both the
string-averaging projections (SAP) and the block-iterative projections (BIP)
methods with fixed strings or blocks, respectively, and prove convergence of
the string-averaging method in the inconsistent case by translating it into a
fully sequential algorithm in the product space.

They consider the successive projections iterative process%
\begin{equation}
x^{0}\in V\text{ is an arbitrary starting point},\text{ }x^{k+1}%
=P_{k(\operatorname{mod}m)+1}\left(  x^{k}\right)  \text{ for all }k\geq0,
\label{censor2003 disp 1}%
\end{equation}
and offer an extension of Gubin%
\index{Gubin}
Polyak%
\index{Polyak}
and Raik's%
\index{Raik}
Theorem \ref{Gubin1967 mark 1} above by replacing the demand that one of the
sets of the CFP is bounded by a weaker condition, as follows ($V$ stands for
the Euclidean space).

\begin{theorem}
\label{censor2003 mark 2}\cite[Theorem 4.4]{censor2003} Let $C_{1}%
,C_{2},\cdots,C_{m}$ be nonempty closed convex subsets of $V.$ If for at least
one set (for explicitness, say $C_{1}$) the cyclic subsequence (of points in
$C_{1}$) $\left\{  x^{km+1}\right\}  _{k\geq0}$ of a sequence $\left\{
x^{k}\right\}  _{k\geq0}$, generated by (\ref{censor2003 disp 1}), is bounded
for at least one $x^{0}\in%
%TCIMACRO{\U{211d} }%
%BeginExpansion
\mathbb{R}
%EndExpansion
^{n}$ then there exist points $x^{\ast,i}\in C_{i},$ $i=1,2,\cdots,m,$ such
that $P_{i+1}\left(  x^{\ast,i}\right)  =x^{\ast,i+1},$ $i=1,2,\cdots,\left(
m-1\right)  ,$ and $P_{1}\left(  x^{\ast,m}\right)  =x^{\ast,1},$ and for
$i=1,2,\cdots,m,$ we have%
\begin{align}
\underset{k\rightarrow\infty}{\lim}x^{km+i+1}-x^{km+i}  &  =x^{\ast
,i+1}-x^{\ast,i},\label{censor2003 disp 2}\\
\underset{k\rightarrow\infty}{\lim}x^{km+i}  &  =x^{\ast,i},
\label{censor2003 disp 3}%
\end{align}
where $\left\{  x^{k}\right\}  _{k\geq0}$ is any sequence generated by
(\ref{censor2003 disp 1}).
\end{theorem}

The SAP method has its origins in \cite{ceh01}, for a recent work about it
see, e.g., Reich%
\index{Reich}
and Zalas%
\index{Zalas}
\cite{Reich2016}.

As in Censor%
\index{Censor}
and Tom%
\index{Tom}
\cite[Page 545]{censor2003}), each string $I_{t}$ is a finite nonempty subset
of $\left\{  1,2,\cdots,m\right\}  ,$ for $t=1,2,\cdots,S,$ of the form
$I_{t}=\left(  i_{1}^{t},i_{2}^{t},\cdots,i_{\gamma\left(  I_{t}\right)  }%
^{t}\right)  ,$ where the length of the string $I_{t},$ denoted by
$\gamma\left(  I_{t}\right)  $, is the number of elements in $I_{t}.$ The
projection along the string $I_{t}$ operator is defined as the composition of
projections onto the sets indexed by $I_{t},$ that is, $T_{t}:=P_{i_{\gamma
\left(  I_{t}\right)  }^{t}}\cdots P_{i_{2}^{t}}P_{i_{1}^{t}}$ for
$t=1,2,\cdots,S.$ Given a positive weight vector $\omega\in%
%TCIMACRO{\U{211d} }%
%BeginExpansion
\mathbb{R}
%EndExpansion
^{S},$ i.e., $\omega_{t}>0,$ $t=1,2,...,S,$ and $\sum_{t=1}^{S}\omega_{t}%
=1,$define the algorithmic operator%
\begin{equation}
T=\sum_{t=1}^{S}\omega_{t}T_{t}, \label{eq:SAP}%
\end{equation}
yielding the SAP method that employs the iterative process%
\begin{equation}
x^{0}\in V\text{ is an arbitrary starting point},\text{ }x^{k+1}=T\left(
x^{k}\right)  \text{ for all }k\geq0. \label{alg:SAP}%
\end{equation}
The convergence of the string-averaging method in the possibly inconsistent
case is included in the following theorem.

\begin{theorem}
\label{censor2003 mark 4}\cite[Theorem 5.2]{censor2003} Let $C_{1}%
,C_{2},\cdots,C_{m},$ be nonempty closed convex subsets of $V$. If for at
least one $x^{0}\in V$ the sequence $\left\{  x^{k}\right\}  _{k\geq0},$
generated by the string-averaging algorithm (Algorithm \ref{alg:SAP}) with $T$
as in (\ref{eq:SAP})), is bounded then it converges for any $x^{0}\in V.$
\end{theorem}

\subsection{2004: Steered sequential projections for the inconsistent convex
feasibility problem}

Censor%
\index{Censor}%
, De Pierro%
\index{De Pierro}
and Zaknoon%
\index{Zaknoon}
\cite{Censor2004} study a steered sequential gradient algorithm which
minimizes the sum of convex functions by proceeding cyclically in the
directions of the negative gradients of the functions and using steered
step-sizes. They apply this algorithm to the convex feasibility problem by
minimizing a proximity function which measures the sum of the Bregman%
\index{Bregman}
distances to the members of the family of convex sets. The resulting algorithm
is a new steered sequential Bregman%
\index{Bregman}
projection method which generates sequences that converge, if they are
bounded, regardless of whether the convex feasibility problem is or is not
consistent (i.e., feasible). For orthogonal projections and affine sets the
boundedness condition is always fulfilled.

The steering parameters in the algorithm form a sequence $\{\sigma
_{k}\}_{k\geq0}$ of real positive numbers that must have the following
properties: $\lim_{k\rightarrow\infty}\sigma_{k}=0,$ $\lim_{k\rightarrow
\infty}(\sigma_{k+1}/\sigma_{k})=1,$ and $\sum_{k=0}^{\infty}\sigma
_{k}=+\infty.$ If instead of $\lim_{k\rightarrow\infty}(\sigma_{k+1}%
/\sigma_{k})=1$ one uses $\lim_{k\rightarrow\infty}\sigma_{km+j}/\sigma
_{km}=1,$ for all $1\leq j\leq m-1,$ then the parameters are called
$\mathit{m}$-steering parameters. For minimization of a function $g\left(
x\right)  :=\sum_{i=0}^{m-1}g_{i}\left(  x\right)  $ where $\{g_{i}%
\}_{i=0}^{m-1}$ is a family of convex functions from $%
%TCIMACRO{\U{211d} }%
%BeginExpansion
\mathbb{R}
%EndExpansion
^{n}$ into $%
%TCIMACRO{\U{211d} }%
%BeginExpansion
\mathbb{R}
%EndExpansion
$ which have continuous derivatives everywhere the cyclic gradient method is
as follows

\begin{algorithm}
\label{alg.cyclic.grad}\cite[Algorithm 5]{Censor2004}\ (\textbf{The m-steered
cyclic gradient method).}

\textbf{Initialization: }$x^{0}\in%
%TCIMACRO{\U{211d} }%
%BeginExpansion
\mathbb{R}
%EndExpansion
^{n}$ is arbitrary.

\textbf{Iterative Step: }Given $x^{k}$ calculate the next iterate $x^{k+1}$ by%
\begin{equation}
x^{k+1}=x^{k}-\sigma_{k}\nabla g_{i(k)}(x^{k})\text{.} \label{iterative step}%
\end{equation}

\textbf{Control Sequence:} $\{i(k)\}_{k\geq0}$ is a cyclic control sequence,
i.e., $i(k)$ $=k\operatorname*{mod}m.$

\textbf{Steering Parameters:} The sequence $\{\sigma_{k}\}_{k\geq0}$ is $m$-steering.
\end{algorithm}

The following convergence result holds.

\begin{theorem}
\label{th.conv.cyclic}\cite[Theorem 6]{Censor2004} Let $\{g_{i}\}_{i=0}^{m-1}
$ be a family of functions $g_{i}:%
%TCIMACRO{\U{211d} }%
%BeginExpansion
\mathbb{R}
%EndExpansion
^{n}\rightarrow%
%TCIMACRO{\U{211d} }%
%BeginExpansion
\mathbb{R}
%EndExpansion
$ which are convex and continuously differentiable everywhere, let $g\left(
x\right)  :=\sum_{i=0}^{m-1}g_{i}\left(  x\right)  $ and assume that $g$ has
an unconstrained minimum. If $\left\{  x^{k}\right\}  _{k\geq0}$ is a
bounded\texttt{\ }sequence$,$ generated by Algorithm \ref{alg.cyclic.grad},
then the sequence $\left\{  g\left(  x^{k}\right)  \right\}  _{k\geq0}%
$\ converges to the minimum of $g.$ If, in addition, $g$ has a unique
minimizer then the sequence\ $\left\{  x^{k}\right\}  _{k\geq0}$ converges to
this minimizer.
\end{theorem}

Applying Algorithm \ref{alg.cyclic.grad} with Bregman%
\index{Bregman}
distance functions and using Theorem \ref{th.conv.cyclic} yields a convergence
result for sequential Bregman%
\index{Bregman}
projections onto convex sets in the inconsistent case.

\begin{algorithm}
\label{alg.cyclic.breg}\cite[Algorithm 14]{Censor2004} \textbf{(Steered cyclic
Bregman%
\index{Bregman}
projections).}

\textbf{Initialization: }$x^{0}\in S$ is arbitrary.

\textbf{Iterative Step: }Given $x^{k}$ calculate the next iterate $x^{k+1}$ by%
\begin{equation}
x^{k+1}=x^{k}+\sigma_{k}\nabla^{2}f(x^{k})\left(  P_{Q_{i(k)}}^{f}%
(x^{k})-x^{k}\right)  \text{.} \label{itierative lease square approximation}%
\end{equation}

\textbf{Control Sequence:} $\{i(k)\}_{k\geq0}$ is a cyclic control sequence,
i.e., $i(k)$ $=k\operatorname*{mod}m.$

\textbf{Steering Parameters:} The sequence $\{\sigma_{k}\}_{k\geq0}$ is an $m
$-steering sequence.
\end{algorithm}

For the Bregman%
\index{Bregman}
function $f\left(  x\right)  :=\left(  1/2\right)  \Vert x\Vert^{2},$ the
algorithm's iterative process takes the form%
\begin{equation}
x^{k+1}=x^{k}+\sigma_{k}\left(  P_{i(k)}(x^{k})-x^{k}\right)  \text{.}%
\end{equation}

Another Bregman%
\index{Bregman}
function is $f\left(  x\right)  =-\operatorname*{ent}x,$ where
$\operatorname*{ent}x$ is Shannon's entropy function which maps the
nonnegative orthant $%
%TCIMACRO{\U{211d} }%
%BeginExpansion
\mathbb{R}
%EndExpansion
_{+}^{n}$ into $%
%TCIMACRO{\U{211d} }%
%BeginExpansion
\mathbb{R}
%EndExpansion
$ by $\operatorname*{ent}x:=-\sum_{j=1}^{n}x_{j}\log x_{j},$ where
\textquotedblleft$\log$\textquotedblright\ denotes the natural logarithms and,
by definition, $0\log0=0.$ The steered cyclic entropy projections method that
is obtained uses the iterative process%
\begin{equation}
x^{k+1}=x^{k}+\sigma_{k}\left(
\begin{array}
[c]{cccc}%
\frac{\displaystyle1}{\displaystyle x_{1}} & 0 & \cdots & 0\\
0 & \frac{\displaystyle1}{\displaystyle x_{2}} & \cdots & 0\\
\vdots & \vdots &  & \vdots\\
0 & 0 & \cdots & \frac{\displaystyle1}{\displaystyle x_{n}}%
\end{array}
\right)  \left(  P_{i(k)}^{f}(x^{k})-x^{k}\right)  \label{entrop step}%
\end{equation}
and its convergence along the lines described above is obtained.

\subsection{2006: Alternating Bregman proximity operators}

It all began with Bregman's paper%
\index{Bregman}
\cite{Bregman1967} which was actually Lev Bregman's Ph.D. work. This paper had
no follow up in the literature until 14 years later in Censor%
\index{Censor}
and Lent%
\index{Lent}
\cite{Censor1981}. For a bibliographic brief review on Bregman functions,
distances and projections consult page 1233 (in the notes and references
section) of the book by Facchinei%
\index{Facchinei}
and Pang%
\index{Pang}
\cite{FacchineiAndPang2003}. See Subsection \ref{subsec:prox2001} here for
additional pointers and read the recent excellent report of Reem%
\index{Reem}%
, Reich%
\index{Reich}
and De Pierro%
\index{De Pierro}
\cite{Reem2018}.

In an attempt to apply compositions of Bregman projections to two disjoint
convex sets Bauschke%
\index{Bauschke}%
, Combettes%
\index{Combettes}
and Noll%
\index{Noll}
\cite{Bauschke2006} investigated the proximity properties of Bregman%
\index{Bregman}%
\ distances. This investigation lead to the introduction of a new type of
proximity operator which complements the usual Bregman proximity operator.

The lack of symmetry inherent to the Bregman distance $D\left(  x,y\right)  $
prompted the authors to consider two single-valued operators defined on $U$,
namely,%
\begin{equation}%
\begin{array}
[c]{l}%
\overleftarrow{prox}_{\varphi}:y\rightarrow\underset{x\in U}{\arg\min}\text{
}\varphi\left(  x\right)  +D\left(  x,y\right)  ,\\
\overrightarrow{prox}_{\psi}:x\rightarrow\underset{y\in U}{\arg\min}\text{
}\psi\left(  y\right)  +D\left(  x,y\right)  ,
\end{array}
\label{Bauschke2006Disp12}%
\end{equation}
where a variaty of assumptions (consult Bauschke%
\index{Bauschke}%
, Combettes%
\index{Combettes}
and Noll%
\index{Noll}
\cite{Bauschke2006}) apply to the functions $\varphi\left(  x\right)  ,$
$\psi\left(  y\right)  $ and $D\left(  x,y\right)  .$ They proposed the
following iterative process \cite[Equation (13)]{Bauschke2006}%
\begin{equation}
\left\{
\begin{array}
[c]{l}%
\text{fix }x_{0}\in U\text{ and set}\\
(\forall n\in%
%TCIMACRO{\U{2115} }%
%BeginExpansion
\mathbb{N}
%EndExpansion
)\text{ \ }y_{n}=\overrightarrow{prox}_{\psi}\left(  x_{n}\right)  \text{ and
}x_{n+1}=\overleftarrow{prox}_{\varphi}\left(  y_{n}\right)  ,
\end{array}
\right.  \label{Bauschke2006Disp 13 and 46}%
\end{equation}
for which they proved convergence that yielded the following corollary.

\begin{corollary}
\cite[Corollary 4.7]{Bauschke2006} Let $A$ and $B$ be closed convex sets in $%
%TCIMACRO{\U{211d} }%
%BeginExpansion
\mathbb{R}
%EndExpansion
^{J}$ such that $A\cap U\neq\varnothing$ and $B\cap U\neq\varnothing$. Suppose
that the solution set $S$ of the problem%
\begin{equation}
\text{minimize }D\text{ over }(A\times B)\cap(U\times U)
\end{equation}
is nonempty. Then the sequence $\left(  \left(  x_{n},y_{n}\right)  \right)
_{n\in%
%TCIMACRO{\U{2115} }%
%BeginExpansion
\mathbb{N}
%EndExpansion
}$ generated by the alternating left-right projections algorithm%
\begin{equation}
\left\{
\begin{array}
[c]{l}%
\text{fix }x_{0}\in U\text{ and set}\\
(\forall n\in%
%TCIMACRO{\U{2115} }%
%BeginExpansion
\mathbb{N}
%EndExpansion
)\text{ \ }y_{n}=\overrightarrow{P}_{B}\left(  x_{n}\right)  \text{ and
}x_{n+1}=\overleftarrow{P}_{A}\left(  y_{n}\right)  ,
\end{array}
\right.
\end{equation}
converges to a point in $S$.
\end{corollary}

\subsection{2012: There is no variational characterization of the cycles in
the method of periodic projections\label{subsec:no-variational}}

Baillon%
\index{Baillon}%
, Combettes%
\index{Combettes}
and Cominetti%
\index{Cominetti}
\cite{Baillon2012} studied the behavior of the sequences generated by periodic
projections onto $m\geq3$ closed convex subsets of a Hilbert space
$\mathcal{H}$. For an ordered family of nonempty closed convex subsets
$C_{1},C_{2},\cdots,C_{m}$ of $\mathcal{H}$ with associated projection
operators $P_{1},P_{2},\cdots,P_{m}$ consider sequences defined by the
following rule: Choose any $x_{0}\in\mathcal{H}$. For every $n=0,1,2,3,\cdots$
perform
\begin{equation}
\left\{
\begin{array}
[c]{l}%
x_{mn+1}=P_{m}x_{mn}\\
x_{mn+2}=P_{m-1}x_{mn+1}\\
\vdots\\
x_{mn+m}=P_{1}x_{mn+m-1}.
\end{array}
\right.  \label{Baillon2012 disp 1}%
\end{equation}

They define the set of cycles associated with the given $m$ closed convex
subsets by:%

\begin{equation}
cyc\left(  C_{1},C_{2},\cdots,C_{m}\right)  =\left\{
\begin{array}
[c]{c}%
\left(  \overline{y}_{1},\overline{y}_{2},\cdots,\overline{y}_{m}\right)
\in\mathcal{H}^{m}\ \text{\ such that}\\
\overline{y}_{1}=P_{1}\overline{y}_{2},\cdots,\ \overline{y}_{m-1}%
=P_{m-1}\overline{y}_{m},\ \overline{y}_{m}=P_{m}\overline{y}_{1}%
\end{array}
\right\}
\end{equation}
and asked if there exist a function $\Phi:\mathcal{H}^{m}\rightarrow%
%TCIMACRO{\U{211d} }%
%BeginExpansion
\mathbb{R}
%EndExpansion
$ such that, for every ordered family of nonempty closed convex subsets
$\left(  C_{1},C_{2},\cdots,C_{m}\right)  $ of $\mathcal{H}$, $cyc\left(
C_{1},C_{2},\cdots,C_{m}\right)  $ can be characterized as the solution set of
a minimization problem of $\Phi?$ They proved the negative answer to this
question by the following theorem.

\begin{theorem}
\label{Baillon2012 mark 2}\cite[Theorem 2.3]{Baillon2012} Suppose that
$\dim\mathcal{H}\geq2$ and let $m$ be an integer at least equal to $3$. There
exists no function $\Phi:\mathcal{H}^{m}\rightarrow%
%TCIMACRO{\U{211d} }%
%BeginExpansion
\mathbb{R}
%EndExpansion
$ such that, for every ordered family of nonempty closed convex subsets
$\left(  C_{1},C_{2},\cdots,C_{m}\right)  $ of $\mathcal{H}$, $cyc\left(
C_{1},C_{2},\cdots,C_{m}\right)  $ is the set of solutions to the variational
problem: minimize$\{\Phi(y_{1},y_{2},...,y_{m})\mid y_{1}\in C_{1},$ $y_{2}\in
C_{2},$ $\cdots,y_{m}\in C_{m}\}.$
\end{theorem}

This shows, in particular, that the Cheney%
\index{Cheney}
and Goldstein%
\index{Goldstein}
\cite{Cheney1959} (see Subsection \ref{subsec:CG} above) result of minimizing
the distance between two disjoint sets cannot be extended to more than two sets.

\subsection{2012: Alternating projections onto two sets with empty or nonempty
intersection\label{subsec:KR}}

Kopeck\'{a}%
\index{Kopeck\'{a}}
and Reich%
\index{Reich}
\cite{Kopecka2012} used tools from nonexpansive operators theory to prove that
alternating projections onto two closed and convex subsets of a real Hilbert
space $H$, generate two subsequences such that the sequence of distances
between them converges to the distance between the two sets. Formally, let
$P_{1}:H\rightarrow S_{1}$ and $P_{2}:H\rightarrow S_{2}$ be the orthogonal
projections of $H$ onto $S_{1}$ and $S_{2},$\ respectively, and denote the
distance between them by $d\left(  S_{1},S_{2}\right)  .$ Define the sequence
$\left\{  x_{n}\mid n\in%
%TCIMACRO{\U{2115} }%
%BeginExpansion
\mathbb{N}
%EndExpansion
\right\}  $ by%
\begin{equation}
x_{2n+1}=P_{1}x_{2n}\text{ and }x_{2n+2}=P_{2}x_{2n+1},
\label{Kopecka2012 disp 1}%
\end{equation}
then the following theorem holds

\begin{theorem}
\label{Kopecka2012 mark 1}\cite[Theorem 1.4]{Kopecka2012}\ Let $S_{1}$ and
$S_{2}$ be two nonempty, closed and convex subsets of a real Hilbert space
$\left(  H,\left\langle \cdot,\cdot\right\rangle \right)  ,$ with induced norm
$\left\Vert \cdot\right\Vert ,$ and let $P_{1}:H\rightarrow S_{1}$ and
$P_{2}:H\rightarrow S_{2}$ be the corresponding nearest point projections of
$H$ onto $S_{1}$ and $S_{2},$\ respectively. Let the sequence $\left\{
x_{n}\mid n\in%
%TCIMACRO{\U{2115} }%
%BeginExpansion
\mathbb{N}
%EndExpansion
\right\}  $ be defined by (\ref{Kopecka2012 disp 1}) then%
\begin{equation}
\underset{n\rightarrow\infty}{\lim}\left\Vert x_{2n+2}-x_{2n+1}\right\Vert
=\underset{n\rightarrow\infty}{\lim}\left\Vert x_{2n+1}-x_{2n}\right\Vert
=d\left(  S_{1},S_{2}\right)  . \label{Kopecka2012 disp 2}%
\end{equation}

\end{theorem}

\subsection{2016: Inconsistent split feasibility problems}

The split common fixed point problem (SCFPP), first proposed in Censor%
\index{Censor}
and Segal%
\index{Segal}
\cite{Censor2009}, requires to find a common fixed point of a family of
operators in one space such that its image under a linear transformation is a
common fixed point of another family of operators in the image space. This
generalizes the convex feasibility problem (CFP), the two-sets split
feasibility problem (SFP) and the multiple sets split feasibility problem (MSSFP).

\begin{problem}
\label{prob:mscfpp}\textbf{The split common fixed point problem.}\newline
Given operators $U_{i}:R^{N}\rightarrow R^{N}$, $i=1,2,\ldots,p,$ and
$T_{j}:R^{M}\rightarrow R^{M},$ $j=1,2,\ldots,r,$ with nonempty fixed points
sets $C_{i},$ $i=1,2,\ldots,p$ and $Q_{j},$ $j=1,2,\ldots,r,$ respectively.
The \texttt{split common fixed point problem (SCFPP)} is%
\begin{equation}
\text{find a vector }x^{\ast}\in C:=\cap_{i=1}^{p}C_{i}\text{ such that
}Ax^{\ast}\in Q:=\cap_{i=1}^{r}Q_{j}.
\end{equation}

\end{problem}

Such problems arise in the field of intensity-modulated radiation therapy
(IMRT) when one attempts to describe physical dose constraints and equivalent
uniform dose (EUD) constraints within a single model, see Censor%
\index{Censor}%
, Bortfeld%
\index{Bortfeld}%
, Martin%
\index{Martin}
and Trofimov%
\index{Trofimov}
\cite{Censor2006}. The problem with only a single pair of sets $C$ in $R^{N}$
and $Q$ in $R^{M}$ was first introduced by Censor%
\index{Censor}
and Elfving%
\index{Elfving}
\cite{Censor1994} and was called the split feasibility problem (SFP). They
used their simultaneous multiprojections algorithm (see also Censor and Zenios%
\index{Zenios}
\cite[Subsection 5.9.2]{censor1997Book}) to obtain iterative algorithms to
solve the SFP.

Iiduka%
\index{Iiduka}
\cite{Iiduka2016} discusses the multiple-set split feasibility problem (MSFP)%
\begin{equation}
\text{Find }x^{\ast}\in C:=\cap_{i\in\mathcal{I}}C^{\left(  i\right)  }\text{
such that }Ax^{\ast}\in Q:=\cap_{j\in\mathcal{J}}Q^{\left(  j\right)  },
\end{equation}
where $C^{\left(  i\right)  }\subseteq%
%TCIMACRO{\U{211d} }%
%BeginExpansion
\mathbb{R}
%EndExpansion
^{N}$ for all $i\in\mathcal{I}:=\left\{  1,2,\cdots,I\right\}  $ and
$Q^{\left(  j\right)  }\subseteq%
%TCIMACRO{\U{211d} }%
%BeginExpansion
\mathbb{R}
%EndExpansion
^{M}$ for all $j\in\mathcal{J}:=\left\{  1,2,\cdots,J\right\}  $ are nonempty,
closed and convex, and $A\in%
%TCIMACRO{\U{211d} }%
%BeginExpansion
\mathbb{R}
%EndExpansion
^{M\times N}$ is a matrix$.$ The author \cite[Page 187]{Iiduka2016} introduces
an inconsistent split feasibility problem (IMSFP) in which the CFP in $%
%TCIMACRO{\U{211d} }%
%BeginExpansion
\mathbb{R}
%EndExpansion
^{N}$ is \textquotedblleft hard\textquotedblright\ (called there
\textquotedblleft absolute\textquotedblright) and the CFP in $%
%TCIMACRO{\U{211d} }%
%BeginExpansion
\mathbb{R}
%EndExpansion
^{M}$ is \textquotedblleft soft\textquotedblright\ (called there
\textquotedblleft subsidiary\textquotedblright) as meant in Subsection
\ref{subsec:hard-soft} here. The infeasibility imposed by assuming that
$\left(  \cap_{i\in\mathcal{I}}C^{\left(  i\right)  }\right)  \cap\left(
\cap_{j\in\mathcal{J}}D^{\left(  j\right)  }\right)  =\emptyset$ where
$D^{\left(  j\right)  }:=\left\{  x\in%
%TCIMACRO{\U{211d} }%
%BeginExpansion
\mathbb{R}
%EndExpansion
^{N}\mid Ax\in Q^{\left(  j\right)  }\right\}  .$

For user-chosen weights $\left(  w^{\left(  j\right)  }\right)  _{j\in
\mathcal{J}}\subset\left(  0,1\right)  $ satisfying $\sum_{j\in\mathcal{J}%
}w^{\left(  j\right)  }=1$ Iiduka%
\index{Iiduka}
employs, for all $x\in%
%TCIMACRO{\U{211d} }%
%BeginExpansion
\mathbb{R}
%EndExpansion
^{N},$ the proximity function%
\begin{equation}
f_{D}\left(  x\right)  :=\frac{1}{2}%
%TCIMACRO{\tsum \limits_{j\in\mathcal{J}}}%
%BeginExpansion
{\textstyle\sum\limits_{j\in\mathcal{J}}}
%EndExpansion
w^{\left(  j\right)  }\left\Vert P_{Q^{\left(  j\right)  }}\left(  Ax\right)
-Ax\right\Vert ^{2}, \label{Iiduka 2016 disp 1.2}%
\end{equation}
where $P_{Q^{\left(  j\right)  }}$ is the metric projection onto $Q^{\left(
j\right)  },$ and represents the IMSFP as the constrained minimization problem%
\begin{equation}
\text{Find }x^{\ast}\text{ such that }f_{D}\left(  x^{\ast}\right)
=\underset{}{\min\left\{  f_{D}\left(  x\right)  \mid x^{\ast}\in
C:=\cap_{i\in\mathcal{I}}C^{\left(  i\right)  }\right\}  }. \label{eq:iiduka}%
\end{equation}

Actually, Iiduka%
\index{Iiduka}
uses nonexpansive mappings $T^{\left(  i\right)  }:%
%TCIMACRO{\U{211d} }%
%BeginExpansion
\mathbb{R}
%EndExpansion
^{N}\rightarrow%
%TCIMACRO{\U{211d} }%
%BeginExpansion
\mathbb{R}
%EndExpansion
^{N}$ and defines $C^{\left(  i\right)  }:=Fix\left(  T^{\left(  i\right)
}\right)  .$ He proposes a sophisticated projections method for solving
(\ref{eq:iiduka}) and formulates conditions under which it converges.

\subsection{2017: Best approximation pairs relative to two closed convex sets}

Given two disjoint closed convex sets,\ say $C$ and $Q,$ a best approximation
pair relative to them is a pair of points, one in each set, attaining the
minimum distance between the sets. Cheney%
\index{Cheney}
and Goldstein%
\index{Goldstein}
\cite{Cheney1959} showed that alternating projections onto the two sets,
starting from an arbitrary point, generate a sequence whose two interlaced
subsequences converge to a best approximation pair. See Subsections
\ref{subsec:CG} and \ref{subsec:KR} here. While Cheney%
\index{Cheney}
and Goldstein%
\index{Goldstein}
considered only orthogonal projections onto the sets, related results, using
the averaged alternating reflections (AAR) method and applying it to not
necessarily convex sets, appeared in Bauschke%
\index{Bauschke}%
, Combettes%
\index{Combettes}
and Luke%
\index{Luke}
\cite{bauschke2004finding} and \cite{Luke2008}. Remaining solely on
theoretical ground, Achiya Dax%
\index{Dax}
\cite{Dax2006} explores the duality relations that characterize least norm
problems. He presents a new Minimum Norm Duality (MND) theorem, that considers
the distance between two convex sets. Roughly speaking, it says that the
shortest distance between the two sets is equal to the maximal
\textquotedblleft separation\textquotedblright\ between the sets, where the
term \textquotedblleft separation\textquotedblright\ refers to the distance
between a pair of parallel hyperplanes that separates the two sets.

The problem of best approximation pair relative to two sets cannot be extended
to more than two sets in view of Baillon%
\index{Baillon}%
, Combettes%
\index{Combettes}
and Cominetti's%
\index{Cominetti}
result \cite{Baillon2012}, see Subsection \ref{subsec:no-variational} here.
From a practical point of view, the best approximation pair relative to two
sets obviously furnishes a solution to the hard and soft constrained
inconsistent feasibility problem, see Subsection \ref{subsec:hard-soft} here.
However, the algorithmic approach is hindered by the need to perform
projections onto the two sets $C$ and $Q$ if they are not \textquotedblleft
simple to project onto\textquotedblright.

Aharoni%
\index{Aharoni}%
, Censor%
\index{Censor}
and Jiang%
\index{Jiang}
\cite{AharoniTechRep2017} propose, for the polyhedral sets case, a process
based on projections onto the half-spaces defining the two polyhedra, which
are more negotiable than projections on the polyhedra themselves. A central
component in their proposed process is the Halpern%
\index{Halpern}%
--Lions%
\index{Lions}%
--Wittmann%
\index{Wittmann}%
--Bauschke%
\index{Bauschke}
(HLWB\footnote{This acronym was dubbed in \cite{CensorComputational2006}.})
algorithm for approaching the projection of a given point onto a convex set.

The HLWB algorithm is applied alternatingly to the two polyhedra. Its
application is divided into sweeps --- in the odd numbered sweeps we project
successively onto half-spaces defining the polyhedron $C$, and in even
numbered sweeps onto half-spaces defining the polyhedron $Q$. A critical point
is that the number of successive projections onto each set's half-spaces
increases from sweep to sweep. The proof of convergence of the algorithm is
rather standard in the case that the best approximation pair is unique. The
non-uniqueness case, however, poses some difficulties and its proof is more involved.

\subsection{2018: Replacing inconsistent sets with set enlargements: ART3,
ARM, Intrepid, Valiant}

Searching for a solution to a system of linear equations is a convex
feasibility problem and has led to many different iterative methods. When the
system of linear equations is inconsistent, due to modeling or measurements
inaccuracies, it has been suggested to replace it by a system of pairs of
opposing linear inequalities creating nonempty hyperslabs. Applying projection
methods to this problem can be done by using any iterative method for linear
inequalities, such as the method of Agmon%
\index{Agmon}
\cite{agmon1954relaxation} and Motzkin%
\index{Motzkin}
and Schoenberg%
\index{Schoenberg}
\cite{motzkin1954relaxation} (AMS). However, in order to improve computational
efficiency, Goffin%
\index{Goffin}
\cite{goffin1971finite} proposed to replace projections onto the hyperslabs by
a strategy of projecting onto the original hyperplane (from which the
hyperslab was created) when the current iterate is \textquotedblleft far
away\textquotedblright\ from the hyperslab, and reflecting into the
hyperslab's boundary when the current iterate is \textquotedblleft close to
the hyperslab\textquotedblright\ while keeping the iterate unchanged if it is
already inside the hyperslab.

In \cite{herman1975relaxation} Herman%
\index{Herman}
suggested to implement Goffin's%
\index{Goffin}
strategy by using an additional enveloping hyperslab in order to determine the
\textquotedblleft far\textquotedblright\ and the \textquotedblleft
close\textquotedblright\ distance of points from the hyperplane, resulting in
his \textquotedblleft Algebraic Reconstruction Technique 3\textquotedblright%
\ (ART3) algorithm. In \cite{censor1985automatic} Censor%
\index{Censor}
also embraced the idea of hyperslabs, and defined an algorithmic operator that
implemented Goffin's%
\index{Goffin}
strategy in a continuous manner, resulting in the Automatic Relaxation Method
(ARM). For applications and additional details see Censor%
\index{Censor}
and Herman%
\index{Herman}
\cite{censor1987some} and \cite{censor1988parallel}.

A fundamental question that remained open since then was whether the
hyperslabs approach to handle linear equations and Goffin's%
\index{Goffin}
principle can be applied to general convex sets. This question was recently
studied by Bauschke%
\index{Bauschke}%
, Iorio%
\index{Iorio}
and Koch%
\index{Koch}
in \cite{bauschke2014method}, see also \cite{bauschke2015projection} and
Bauschke, Koch and Phan%
\index{Phan}
\cite{bauschke2015stadium} for further details and interesting applications.
They defined convex sets enlargements instead of hyperslabs and used them to
generalize the algorithmic operator that appeared in Herman%
\index{Herman}
\cite{herman1975relaxation}. They defined an operator which they called the
\textquotedblleft intrepid projector\textquotedblright, intended to generalize
the ART3 algorithm of \cite{herman1975relaxation} to convex sets. Motivated by
\cite{bauschke2014method}, Censor%
\index{Censor}
and Mansour%
\index{Mansour}
\cite{Censor2018} present a new operator, called the \textquotedblleft valiant
operator\textquotedblright, that enables to implement the algorithmic
principle embodied in the ARM of \cite{censor1985automatic} to general convex
feasibility problems. Both ART3 and ARM seek a feasible point in the
intersections of the hyperslabs and so their generalizations to the convex
case seek feasibility of appropriate enlargement sets that define the extended
problem.\bigskip

\textbf{Acknowledgement}. We thank Heinz Bauschke for calling our attention to
some relevant references. This work was supported by Research Grant No.
2013003 of the United States-Israel Binational Science Foundation (BSF).%

%TCIMACRO{\TeXButton{Place Index Here}{\printindex}}%
%BeginExpansion
\printindex
%EndExpansion

\bibliographystyle{plain}
\bibliography{acompat,Maroun}

\end{document}